\documentclass[a4paper,notitlepage, 12pt,reqno]{article}

  \usepackage{mathtext}
 \usepackage[cp1251]{inputenc}
  \usepackage[T2A]{fontenc}
 \usepackage[russian,english]{babel}
 \usepackage[all,arc,poly,2cell,curve,arrow,tips]{xy}
\usepackage{graphicx}

  \usepackage{enumerate, eucal, amsthm,amsmath, amssymb}

  \allowdisplaybreaks[4]

 \theoremstyle{definition}
 \newtheorem{defn}{Definition}

 \theoremstyle{plain}
 \newtheorem{thm}{Theorem}
 \newtheorem*{thm*}{Theorem}
 \newtheorem{prop}{Proposition}
  \newtheorem*{prop*}{Предложение}
 
  \newtheorem*{cor*}{Следствие}
 
  \newtheorem*{lem*}{Лемма}

 \theoremstyle{remark}
 
 \newtheorem{remark}[defn]{Remark}
 \newtheorem*{remark*}{Замечание}

 \renewcommand{\abstractname}{}

  \newcounter{ab}
\title{ Central elements in the universal enveloping algebra and function
of matrix elements. \footnote
 {The work was supported by grant MK-4594.2013.1.}}

 \author{D.V. Artamonov, V.A. Golubeva}
  \date{}

\begin{document}
 \maketitle

\renewcommand{\abstractname}{}

\begin{abstract}
In the paper a construction of central elements in
$U(\mathfrak{o}_N)$ and $U(\mathfrak{g}_2)$ based on invariant theory
is given.  New function of matrix elements that appear in
description of the center of $U(\mathfrak{g}_2)$ are defined.

\end{abstract}

\section{Introduction}

 Let $\mathfrak{g}$ be a simple Lie algebra  and let  $V$  be it's standard representation,  put  $dim V=N$.
 Central elements in the universal enveloping algebra $U(\mathfrak{g})$
can be expressed as functions of matrix elements of the matrix $L$
that is defined below, for different algebras different functions
(determinants, pfaffians, hafnians) are used \cite{i1}, \cite{i2},
\cite{Molev}. However in cited papers  the mentioned functions of
matrix elements are not derived from some required properties. The
corresponding formulas are presented and then it is proved that
define central elements.

In the present paper a general construction of central elements is
given. In the case of orthogonal algebra it leads to formulas
involving pfaffians. Also this scheme is applied  to the exceptional
algebra $\mathfrak{g}_2$,
 in this case we obtain some new functions of matrix elements.

Let us define the matrix $L$. Let
 $v_{\alpha}$ be a base in $\mathfrak{g}$, and $V_{\alpha}$ be a matrix that correspond to $v_{\alpha}$ in the standard representation.
  Put

 $$L= V_{\alpha}\otimes v_{\alpha}.\footnote{Everywhere in the paper the summation by the repeating indices is suggested}$$

Note that $L$ can be considered as an $N\times N$ matrix, whose
elements
 belong to $\mathfrak{g}$. Then the following facts take place

\begin{enumerate}
\item In the case when $\mathfrak{g}$ belongs to the series  $A$, in
$U(\mathfrak{g})$ there are central elements given by the formula

\begin{equation}
\label{f1}
 C_k=\sum_{|I|=k}Det L_{I},\,\,\,\, k=2,...,N
\end{equation}

where  $I\subset \{1,...,N\}$, and $L_I$ is a submatrix in $L$,
defined as  $L_I=(L_{ij})_{i,j\in I}$, and $Det$ is a double
determinant (see formula \eqref{2}). The summation is taken over all
subsets $I\subset\{1,...,N\}$ that consist of $k$ elements.

In particular the element

\begin{align}
\begin{split}
\label{2}
 Det(L)=\frac{1}{N!}\sum_{\sigma,\sigma'\in
 S_N}L_{\sigma(1),\sigma'(1)}...L_{\sigma(N),\sigma'(N)}
\end{split}
\end{align}

 is central.



 \item  In the case when $\mathfrak{g}$ belongs to the series $B$ or $D$,
in $U(\mathfrak{g})$ there are the following matrix elements

\begin{equation}
C_2k=\sum_{|I|=k}(Pf L_I)^2,\,\,\,\,
\end{equation}

 where  $k=1,...,N/2$  in the case of series $D$ and $k=1,...,(N-1)/2$  in the case of series
 $B$.

Also in the case $D$ the element

\begin{align}
\begin{split}
PfL=\frac{1}{N!2^N}\sum_{\sigma\in
S_{2N}}L_{\sigma(1),\sigma(2)}...L_{\sigma(2N-1),\sigma(2N)}
\end{split}
\end{align}

is  central.

\item In the case when $\mathfrak{g}$ belongs to the series $C$
in the universal enveloping algebra there are central elements that
are expressed through the so called hafnians of submatrices of the
matrix $L$ \cite{Molev}.
\end{enumerate}

Thus one can say that in the case of the series $A$ the central
elements are expressed through determinants, in the case of the
series $B$, $D$ the central elements are expressed through
pfaffians. But for exceptional Lie algebras a relation between
central elements and new functions of matrix elements is not pointed
out.


There appears a question. Why in construction of central elements in
the case of series $B$, $D$   the pfaffians and not other function
of matrix elements appear?  Which functions appear in the
construction of central elements in the case of the algebra
$\mathfrak{g}_2$?

\section{The content of the paper}

In the present paper a construction of central elements in the
universal enveloping algebra based on the first main theorem of the
invariant theory is given. For the construction of the central
elements a notion of an $m$-invariant is introduced. An
$m$-invariant is polynomial in variables $m_{i,j}$,
$i,j\in\{1,...,N\}$,  that is invariant under the action of the Lie
algebra $\mathfrak{g}$. The action of the generator $v$ of the Lie
algebra,  to which  in the standard representation there corresponds
the matrix $V$, on these variables is given by formula

$$m_{i,j}\mapsto V_{i,c}m_{c,j}+m_{i,d}V_{j,d}.$$

It is proved in the present paper that in the case of series $A$,
$B$, $C$, $D$, and also in the case $\mathfrak{g}_2$ when one
substitute into an $m$-invariant elements $L_{i,j}$ instead of
$m_{i,j}$, and $L_{i,j}$ are multiplied using the symmetrized
product, one gets a central element in $U(\mathfrak{g})$.

In Section \ref{or}  the cases of series $B$ and $D$ are considered.
Using the first main theorem of the invariant theory a description
of a general $m$-invariant is given.  Then we present new relations
that appear when one substitutes  $L_{i,j}$ instead of $m_{i,j}$. As
a corollary one obtains a well known description of the center of
$U(\mathfrak{o}_N)$.  Let us stress that in this approach  the
pfaffians in formulas appear very natural from the first main
theorem of the invariant theory.

Secondly in Section \ref{gg} the case of the exceptional Lie algebra
$\mathfrak{g}_2$ is considered.\footnote{ Mention that in paper
\cite{AG1} an expression for central element as sums of squares of
pfaffians was obtained as in the case of orthogonal algebra. } The
central elements for $\mathfrak{g}_2$ are constructed using the
first main theorem of the invariant theory.

In Section \ref{nf} we define some new functions of matrix  elements
through which the central elements are expressed. Let us write some
of these functions. Let $M=(m_{ij})$ be a $8\times 8$  matrix(or
$7\times 7$ matrix),  whose rows and columns are indexed by
octonions (imaginary octonions). Let  $\omega_{i_1,i_2,i_3}$  be
structure constants of octonions. Put

\begin{equation}
 \omega_{i_1,...,i_k}:=skew_{i_1,...,i_k}\omega_{i_1,i_2,s_1}\omega_{s_1,i_3,s_2}...\omega_{s_{k-1},i_k,1},
\end{equation}

where $skew_{i_1,...,i_k}$ denotes antisymmetrization of the indices
$i_1,...,i_k$. This tensor is skewsymmetric. Then the new function
of matrix elements are

\begin{align}
\begin{split}
\label{new}
 G_k^k(L)=\omega_{i_1,...,i_k}\omega_{j_1,...,j_k}
m_{i_1,j_1}...m_{i_{k},j_{k}}
\end{split}
\end{align}

Also in Section \ref{nf} other functions appear.

However when one substitutes $L_{i,j}$ instead of  $m_{i,j}$ into
these functions during the construction of central elements it turns
out to be possible to express functions \eqref{new} through
determinants of submatrices of the matrix $L$.

\section{Preliminaries}

\subsection{Invariant polynomials. The first main theorem of the invariant theory. $m$-invariants.}

Let us be  given a Lie algebra $\mathfrak{g}$, let $V$ be it's
standard representation. Let  $x_1,...,x_m$ be vectors from $V$, put
$dim V=N$, and denote as $x_k^i$, $i=1,...,N$ the coordinates of
vectors $x_k$.

The first main theorem of the invariant theory is a theorem that
describes  generators in the algebra of  polynomials in variables
$x_k^i$, that are invariant under the action of $\mathfrak{g}$.

We call an invariant polynomial in variables $x_k^i$  an
$x$-invariant.

Let us be given an  $N\times N$ matrix   $M$, whose elements are
variables $m_{ij}$ with no relations between them.  Define an action
on the variables  $m_{ij}$ of the element $v$ algebra $\mathfrak{g}$
by formulas

\begin{equation}
\label{f2}
 M\mapsto VM+MV^t, \text{ or } m_{ij}\mapsto
V_{ik}m_{kj}+m_{jl}V_{il},
\end{equation}

where $v\in \mathfrak{g}$,  and $V_{i,j}$ is a matrix, corresponding
to $v$ in the standard representation.

\begin{defn} A polynomial $m_{ij}$, that is invariant under the action of the algebra is called
an  $m$-invariant.
\end{defn}

There exist an obvious correspondence  between  homogeneous
$x$-invariant of degree $2k$ and homogeneous $m$-invariant of degree
$k$. To an $x$-invariant
$$t^{i_1,j_1...,i_k,j_k}x^1_{i_1}x^{1'}_{i'_1}...x^k_{i_k}x^{k'}_{i'_k}$$
there corresponds an $m$-invariant
$$t^{i_1,j_1...,i_k,j_k}m_{i_1,j_1}...m_{i_kj_k}.$$

\section{A relation between $m$-invariants and central elements in the universal enveloping algebra}

Let us for an algebra $\mathfrak{g}$ construct a matrix $L$ whose
elements belong to $\mathfrak{g}$. Let $v_{\alpha}$ be a base in
$\mathfrak{g}$ and denote as $V_{\alpha}$ a matrix corresponding to
$v_{\alpha}$ in the standard representation. Put

$$L=V_{\alpha}\otimes v_{\alpha}.$$

In the case $\mathfrak{g}=\mathfrak{o}_N$ the matrix $L$  up to
multiplication  by a constant equals to
$$(F_{i,j}),\,\,\, i,j=1,...,N,$$  where generators  $F_{i,j}$  are defined by formula $F_{i,j}=E_{i,j}-E_{j,i}.$ Thus the matrix  $L$ is skew-symmetric.

In the case $\mathfrak{g}=\mathfrak{g}_2$ the matrix $L$  up to
multiplication  by a constant equals to $$(G_{i,j}),\,\,\,
i,j=1,...,7,$$ where generators $G_{i,j}$ are defined as follows.
The Lie algebra $\mathfrak{g}_2$ is the algebra of differentiations
of octonions. For octonions $x,y$ define a differentiation $G_{x,y}$
that act on an octonion $z$ as follows

$$G_{x,y}(z)=[[x,y],z]-3[x,y,z]\text{ where }[x,y,z]=(xy)z-x(yz).$$

Take a standard base $1,e_1,...,e_7$ in the algebra of octonions,
such that $e_1,...,e_7$ are standard imaginary octonions, take
$x=e_i, y=e_j$ and put
$$G_{i,j}=G_{e_i,e_j}.$$

In the case of both algebras $\mathfrak{o}_N$ and $\mathfrak{g}_2$
commutation relations between generators can be written in a similar
way. For $g\in\mathfrak{o}_N$ define the element $F_{g\{i_1,i_2\}}$
as follows.
 Identify a pair of indices $\{i_1,i_2\}$ with the wedge-product of
 vectors of standard representation
$e_{i_1}\wedge e_{i_2}$. Let  $$g(e_{i_1}\wedge
e_{i_2})=c_{\alpha}e_{p_{\alpha}}\wedge e_{q_{\alpha}}.$$

Then put
$$F_{g\{i_1,i_2\}}=c_{\alpha}F_{p_{\alpha},q_{\alpha}}.$$

In \cite{AG} it is shown that commutation  relations between
generators of $\mathfrak{o}_N$ can be written as follows

\begin{equation}
\label{cf}
[g,F_{i,j}]=F_{g\{i,j\}},\,\,\, g\in\mathfrak{o}_N
\end{equation}

Analogously  in \cite{AG1} it is shown that in the case
$\mathfrak{g}_2$ one has

\begin{equation}
\label{cg}
[g,G_{i,j}]=G_{g\{i,j\}},\,\,\, g\in\mathfrak{g}_2.
\end{equation}

Below we denote as $L_{i,j}$  the generators $F_{i,j}$ in the case
$\mathfrak{o}_N$ and generators $G_{i,j}$ in the case
$\mathfrak{g}_2$.

Let us prove Proposition.

\begin{prop}In the case of algebras $\mathfrak{g}=\mathfrak{o}_N$, $\mathfrak{g}_2$ if $$T=t^{i_1,j_1...,i_k,j_k}m_{i_1,j_1}...m_{i_kj_k}$$ is an
$m$-invariant, then
$$\overline{T}=t^{i_1,j_1...,i_k,j_k}L_{i_1,j_1}\cdot...\cdot
L_{i_kj_k}$$  is a central element in $U(\mathfrak{g})$. Here
 $\cdot$ is a symmetrized product.
\end{prop}

\proof

There exist a mapping

\begin{equation}
\label{ot}
 m_{ij}\mapsto L_{ij}
\end{equation}

from the algebra of polynomials in variables $m_{ij}$ into the
algebra of polynomials in variables $L_{ij}$ with the symmetrized
product. The last algebra is embedded in $U(\mathfrak{g})$.

According to formulas \eqref{cf}, \eqref{cg} the mapping given by the formula \eqref{ot},
transforms the action of an element $g\in \mathfrak{g}$ on the
polynomials $m_{ij}$ into the operation of commutation $[g,.]$.
Hence the invariant polynomials are mapped into central elements.

\endproof

\section{Orthogonal algebra}
\label{or}

In the Section the first main theorem of the invariant theory is
formulated. Using this theorem generators in the algebra of
$m$-invariants are written. These generators are encoded by graphs.
Relations between these generators are written. These relations
allow to define basic $m$-invariants and to express all invariant
though basic invariants.

\subsection{The first main theorem}

\begin{thm}(see \cite{W})
Let $x_1,...,x_m$be vectors of the standard representation $V$ of
the algebra $\mathfrak{o}_N$. Then the algebra of polynomials in
coordinates of vectors $x_k$ that are invariant under the action of
$\mathfrak{o}_N$ is generated by $(x_k,x_l)=\sum_i x_k^i x_l^i$ for
all $k,l$, and also in the case $N=m$, by the polynomial
$det(x_k^i)=det[x_1,...,x_N]$, where $[x_1,...,x_N]$ is a matrix
that is constructed from the columns of coordinates of vectors
$x_1,...,x_N$.
\end{thm}

\subsection{Examples of $m$-invariants.}

Using the correspondence between $x$ and $m$-invariants let us construct examples of  $m$-invariants.

\label{inv}

\subsubsection{The trace}
\label{i1} Take as an $x$-invariant the scalar product. That is put
$T(x_1,x_2)=x_1^{i}x_2^{i}$. Then the corresponding $m$-invariant is
the trace of $M$.

\subsubsection{The pfaffian}
\label{i2}
 Take as an $x$-invariant the determinant. That is
$T(x_1,...,x_{k})=det[x^1,...,x^k]$. The components of $T$ are the
following
$$t^{i_1,j_1...,i_k,j_k}=0, \text{ if }\{i_1,j_1...,i_k,j_k\}\neq\{1,...,2k\},$$
$$t^{i_1,j_1...,i_k,j_k}=(-1)^{\sigma}, \text{ if }(i_1,j_1...,i_k,j_k)=(\sigma(1),...,\sigma(2k)).$$

The corresponding $m$-invariant equals

$$k!2^kPfM.$$

\subsubsection{The determinant} \label{i3}
 Take the following  $x$-invariant

$$T(x_1,x_{1'},...,x_{k},x_{k'})=det[x_1,...,x_k]det[x_{1'}...x_{k'}].$$
In the components one gets
$$t^{i_1,j_1...,i_k,j_k}=0, \text{ if }\{i_1,...,i_k,\}\neq\{1,...,k\}, \text{ or } \{j_1,...,j_k,\}\neq\{1,...,k\}$$
$$t^{i_1,j_1...,i_k,j_k}=(-1)^{\sigma}(-1)^{\sigma'}, \text{ if
 }(i_1,...,i_k)=(\sigma(1),...,\sigma(2)),\,\,\,\,(j_1...,j_k)=(\sigma'(1),...,\sigma'(2k))$$

The corresponding $m$-invariant equals

$$k!detM.$$

\subsection{The graphical description of $m$-invariants}

Let us describe a general  $m$-invariant. The general $x$-invariant
is a linear combination of traces and determinants. Let us find an
$m$-invariant that correspond to such $x$-invariant.

A general $m$-invariant  is encoded by an oriented graph with some
additional information of the following kind. The vertices are of
colored into two colors: white or black. Every black vertex belongs
to exactly $n$ edges the number $n$  is the same for all edges. An
order on the edges that belong to a black vertex is fixed. A white
vertex belongs to two edges. All edges are numerated by numbers from
$1$ to $K$.

\begin{figure}
\centering\includegraphics[bbllx=0,bblly=0,bburx=100,bbury=100]{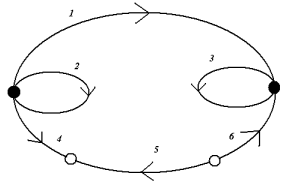}
\caption{A graph that encodes an $m$-invariant} \label{figu1}
\end{figure}

To construct an $m$-invariant we take the product
 $m_{i_1,j_1}...m_{i_K,j_K}$.

To every white vertex there corresponds a contraction.

 In the case when the edges with numbers
$a$, $b$ begin in this vertex one takes the contraction
$$m_{i_a,c}m_{i_b,c}.$$

In the case when the edges with numbers $a$, $b$ end in this vertex
one takes the contraction
$$m_{c,j_a}m_{c,j_b}.$$

\begin{figure}
\label{fig6}
\centering\includegraphics[bbllx=0,bblly=0,bburx=100,bbury=100]{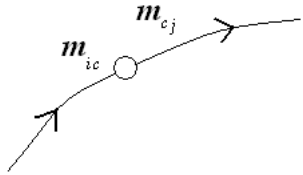}
\caption{A contraction corresponding to a white vertex}
\end{figure}

In the case when the edge $a$ ends in the considered vertex and the
edge with the number $b$ begins in the vertex one takes the
contraction (see figure \ref{fig6})
$$m_{i_a,c}m_{c,j_b}.$$

To every black vertex there corresponds an antisymmetrization.
Suggest that the black vertex belongs to edges with numbers
$a_1$,...,$a_n$. For every edge with the number  $a_{k}$ take the
index $i_{a_k}$ if the edge begins in the vertex and take the index
 $b_{k}$ if the edge ends in the vertex. Then an antisymmetrization of the chosen indices is taken.  The obtained expression is an $m$-invariant,  that corresponds to a graph.

This description of  $m$-invariants follows from the description of
$x$-invariant and the construction of $m$-invariant from
$m$-invariants.

\begin{remark}
An $m$-invariant can be written as  a linear combination of products
of determinants and pfaffians of submatrices of products of $M$ and
$M^t$. However this description is needed below.
\end{remark}

Let us describe graphs that correspond to invariants from the
subsection \ref{inv}.

\begin{figure}
\centering\includegraphics[bbllx=0,bblly=0,bburx=100,bbury=100]{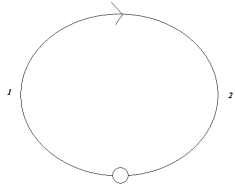}
\caption{A graph that encodes a trace} \label{fig4}
\end{figure}

To the trace there corresponds the graph shown in the figure
\ref{fig4}.

\begin{figure}
\centering\includegraphics[bbllx=0,bblly=0,bburx=100,bbury=100]{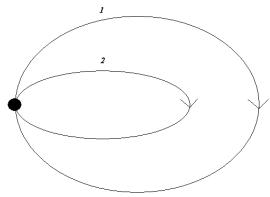}
\caption{A graph that encodes a pfaffian} \label{fig2}
\end{figure}

To the pfaffian in the case $k=2$ there corresponds the graph from
the figure \ref{fig2}

\begin{figure}
\centering\includegraphics[bbllx=0,bblly=0,bburx=100,bbury=100]{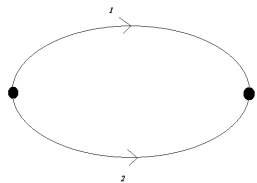}
\caption{A graph that encodes a determinant} \label{fig3}
\end{figure}

To the determinant in the case  $k=2$ there corresponds a graph
shown in the figure \ref{fig3}

\subsection{Relations}

If one puts no relation on the matrix $M=(m_{ij})$ then all
relations between $m$-invariants are corollaries of relations
between $x$-invariant,  these relation are well-known (see the
second main theorem in the invariant theory in \cite{W}).

For the element of the matrix  $L$ there exist the following
relation of skew-symmetry
$$L_{i,j}=-L_{j,i}.$$

Let us put onto the elements $M$  relations of skew symmetry
$m_{i,j}=-m_{j,i}$. In this case in addition to relations following
from the second main theorem of the invariant theory  new relations
between $m$-invariants appear. These relation allow to define basic
invariants and to express other invariant though them.

Let us as formulate these relation.

\subsubsection{The first relation}

In the case of skew symmetric matrix  $M$ of size $2\times 2$ the
following obvious relation takes place
$$(PfM)^2=DetM.$$

\begin{figure}
\label{fig7}
\centering\includegraphics[bbllx=0,bblly=0,bburx=250,bbury=100]{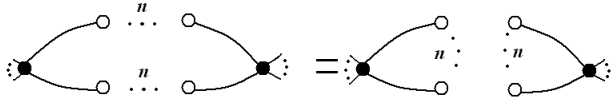}
\caption{The first relation}
\end{figure}

It gives a relation shown in the figure \ref{fig7}. This figure
states the following relation. Let us be given a graph where from
one black vertex to another one go two similar paths though white
vertices (in these paths the same are the numbers of white vertices
and orientations of edges). Then the corresponding $m$-invariant
equals to an $m$-invariant that is defined by the graph where these
paths go not from one black vertex to another but from return to the
same black vertex.

\subsubsection{The second relation}

In the case of matrices  $M$ and $M'$ of size $2\times 2$ one has
\begin{align}
\begin{split}
&Det(M+M')-DetM-DetM'=\frac{1}{2}(m_{1,1}m'_{2,2}-m_{2,1}m'_{1,2}-m_{1,2}m'_{2,1}+m_{2,2}m'_{1,1}+\\
&+m'_{2,2}m_{1,1}-m'_{1,2}m_{2,1}-m'_{2,1}m_{1,2}+m'_{1,1}m_{2,2})
\end{split}
\end{align}

\begin{figure}
\label{fig8}
\centering\includegraphics[bbllx=0,bblly=0,bburx=250,bbury=100]{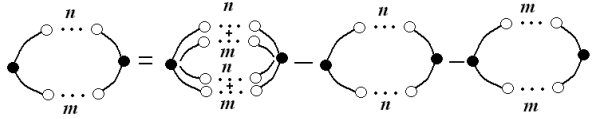}
\caption{The second relation}
\end{figure}

This equality gives a relation shown on the picture \ref{fig8}. On
this picture the first  graph on the right side of the equality
denotes the following invariant. First one takes a $2\times
2$-matrix that corresponds to the path with  $n$ white vertices.
Then  one takes as analogous product  for the path with $m$ white
vertices. Then one takes a sum of these matrices and then one takes
it's determinant.

\subsubsection{The third relation}

In the case of skew-symmetric matrices $M$ and $M'$  of size
$4\times 4$ one has an equality

\begin{align}
\begin{split}
& Pf(M+M')-Pf M -PfM'=\frac{1}{2}(
m_{1,2}m'_{3,4}-m_{1,3}m'_{2,4}+m_{1,4}m'_{2,3}+\\&+m'_{3,4}m_{1,2}-m'_{2,4}m_{1,3}+m'_{2,3}m_{1,4})
\end{split}
\end{align}

\begin{figure}
\label{fig9}
\centering\includegraphics[bbllx=0,bblly=0,bburx=250,bbury=100]{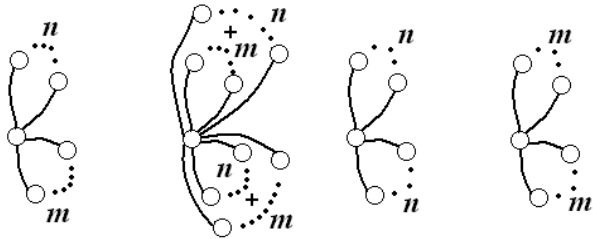}
\caption{The third relation}
\end{figure}

This equality gives us a relation shown in the figure \ref{fig9}.

\subsubsection{The fourth relation}

In the case of skew symmetric matrix  $M$ one has a relation

$$Pf(M^k)=(PfM)^k.$$

\subsubsection{Basic invariants}

Using the  relations above in the case of skew-symmetric matrix $M$
one can express every relation through the $m$-invariant of type

\begin{equation}
\label{bai}
 Tr(M^k)\text{  or }Pf(M^{n_1}+...+M^{n_k}).
\end{equation}

Using standard formulas connecting different symmetric polynomials
one can express the invariants \eqref{bai} through the invariants

\begin{equation}
\sum_{|I|=k}(PfM_I)^2,\text{  or }PfM\,\,\,\, I\subset
\{1,...,N\},\,\,\, M_{I}=(m_{i,j})_{i,j\in I}.
\end{equation}

\section{The algebra $\mathfrak{g}_2$}
\label{gg}

 Let us formulate the first main theorem of the invariant theory in the case of the algebra
  $\mathfrak{g}_2$. Using this theorem let us describe
$m$-invariants and central elements. Some of these central elements
are written explicitly.

\subsection{The first main theorem of the invariant theory}

\begin{thm}(see \cite{S1},\cite{S2})
Let $x_1,...,x_m$ be vectors of the standard representation  $V$ of
the algebra $\mathfrak{g}_2$. Then the algebra of invariant
polynomials in variables $x_1,...,x_m$ is generated by polynomials
$$str(x_{j_1}(x_{j_2}(x_{j_3}...x_{j_r})...))),$$ where  vectors are multiplied as octonions and
 $str$ denotes the symmetrized trace. And the trace denotes the
 operation of taking the real part.
\end{thm}

Let us give another formulation of this theorem. Denote as
$\omega_{i_1,i_2,i_3}$ the structure constants of octonions. Put

\begin{equation}
 \omega_{i_1,...,i_k}:=skew\omega_{i_1,...,i_k}\omega_{i_1,i_2,s_1}\omega_{s_1,i_3,s_2}...\omega_{s_{k-1},i_k,1},
\end{equation}

where $skew_{i_1,...,i_k}$ denotes antisymmetrization of the indices
$i_1,...,i_k$.

One describe this tensor as follows. Identify the index $i_s$ with
the basic octionon $e_{i_s}$, then $\omega_{i_1,...,i_k}$ is a
skew-symmetric $k$-tensor, it's component equals to $0$, if

$$e_{i_1}....e_{i_k}\neq \pm 1,$$

its component equals to $1$, if

$$e_{i_1}....e_{i_k}=1.$$

Using this tensor one can reformulate the first main theorem as
follows

\begin{thm}
Let $x_1,...,x_m$ be vectors of the standard representation $V$ of
the algebra $\mathfrak{g}_2$. Then the algebra of polynomials in
coordinates of vectors $x_1,...,x_m$ is generated by polynomials
$$\omega_{i_1,...,i_k}x_{j_1}^{i_1}....x_{j_k}^{i_k}.$$

\end{thm}

\subsection{$m$-invariant}\label{mg2}

Let us give description of $m$-invariants, using the first main
theorem of the invariant theory.

Consider a product

$$m_{i_1,j_1}...m_{i_k,j_k}.$$

Let us divide the set of indices $\{i_1,j_1,...,i_k,j_k\}$ into $p$
groups  $\{a_1,...,a_s\}$,...,$\{b_1,...,b_t\}$. The indices of the
first group are contracted with $\omega_{a_1,...,a_s}$,...,
 the indices in the last group are contracted with $\omega_{b_1,...,b_t}$.

One can give a graphical description of an $m$-invariant. An
invariant is described by an oriented graph whose edges are
numerated by numbers $1,...,p$. To the edge with the number $l$
there corresponds $m_{i_l,j_l}$. A vertex where the edges with
numbers $a_1,...,a_s$ end end the edges with numbers $b_1,...,b_t$
begin there correspond a contraction
$$j_{a_1},...,j_{a_s},i_{b_1},...,i_{b_t}\text{ с } \omega_{j_{a_1},...,j_{a_s},i_{b_1},...,i_{b_t}}.$$

\begin{figure}
\centering\includegraphics[bbllx=0,bblly=0,bburx=100,bbury=100]{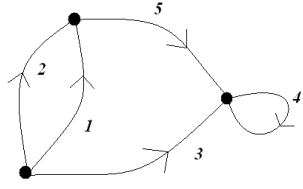}
\caption{A graph that encodes an $\mathfrak{g}_2$-invariant}
\label{fig5}
\end{figure}

\subsection{Central elements}

Let us find, which $m$-invariant give nonzero central elements. In
the case of the algebra $\mathfrak{g}_2$ the elements $L_{ij}$
satisfy the relation \cite{dla}

\begin{equation}
L_{ij}\omega_{i,j,l}=0
\end{equation}

Using this relation and the description of $m$-invariants for
$\mathfrak{g}_2$ from Section \ref{mg2} one finds that the central
elements are describes as follows. One takes a product

$$L_{i_1,j_1}...L_{i_k,j_k},$$

then the set of indices $\{i_1,j_1...,i_k,j_k\}$  is divided into
subsets $\{a_1,...,a_t\}\sqcup...\sqcup\{b_1,...,b_t\}$, where one
set set does not contain $i_l,j_l$ simultaneously.
 Then every set of indices is contracted with the corresponding tensor $\omega$.

\subsection{The functions $G$}
\label{nf}

Let $M=(m_{ij})$ be a  $8\times 8$ (or $7\times 7$) matrix, whose
rows and columns are indexed by octonions (imaginary octionos).

Fix an integer  $k$, $2\leq k \leq 8$,  integers $n_1,...,n_p$,
$n_1+...+n_p=k$  and $m_1,...,m_q$, $m_1+...+m_q=k$

\begin{align}
\begin{split}
\label{new1} &
G_{m_1,...,m_q}^{n_1,...,n_p}(L)=\\&=\sum_{\{i_1,...,i_k\}=I_1\sqcup
... \sqcup I_p,\{j_1,...,j_k\}=J_1\sqcup ... \sqcup
J_q}\omega_{i^1_1,...,i^1_{n_1}} ...\omega_{i^1_p,...,i^p_{n_p}}
\omega_{j^1_1,...,j^1_{m_1}}...\omega_{j^q_1,...,j^q_{m_q}}
L_{i_1,j_1}...L_{i_{k},j_{k}},
\end{split}
\end{align}
where the first summation is taken over all partitions of the sets
$I=\{i_1,...,i_k\}$ и $J=\{j_1,...,j_k\}$ such that
$I_1=\{i_1^1,...,i^1_{n_1}\}$,...,$J_1=\{j_1^1,...,j^1_{m_1}\}$....

As it is proved above  the functions
$G_{m_1,...,m_q}^{n_1,...,n_p}(L)$ generate the algebra of central
elements in $U(\mathfrak{g}_2)$.

However there exist relations between these generators. Actually
when writes the basic central elements the corresponding functions
$G_{m_1,...,m_q}^{n_1,...,n_p}(L)$ can be expressed through the
determinants and pfaffians.

\subsection{Examples}

It is known that in $U(\mathfrak{g}_2)$ the exist primitive central
elements of orders $2$ and $6$.

\subsubsection{The central element of the order $2$.}
Consider the invariant of the order $2$. Is is given by the function
$G_2^2(L)$. To write it let us take the product

\begin{equation}
L_{i_i,j_1}L_{i_2,j_2},
\end{equation}

the indices $i_1,i_2$ must be contracted with $\omega_{i_1,i_2,1}$,
and the indices $j_1,j_2$ must be contracted with
$\omega_{j_1,j_2,1}$.

Since $\omega_{i_1,i_2,i_3}$ are structure constants of octonions,
and squares of base octonions equal to $-1$,  then
$$\omega_{i_1,i_2,1}=-1\Leftrightarrow i_1=i_2.$$

Hence the contraction is actually a sum

\begin{equation}
L_{i_i,j_1}L_{i_1,j_1}.
\end{equation}

Thus the central element is the Casimir element.

\subsubsection{The central element of  higher orders.} Consider the invariant of order
$m$ that is defined by the function $G_m^m(L)$. To write it let us
take the product

\begin{equation}
L_{i_i,j_1}...L_{i_m,j_m},
\end{equation}

indices $i_1,...,i_m$ are contracted with $\omega_{i_1,...,i_m}$, а
indices $j_1,...,j_m$ are contracted with $\omega_{j_1,...,j_m}$.

There exist $7$ imaginary octonions and there product equals $\pm
1$.

 the contraction with $\omega_{i_1,...,i_6}$ can be done as
follows. Thus one obtains that $\omega_{j_1,...,j_m}=0$ for $m=5,6$,
and for $m=7$ the contraction with $\omega_{j_1,...,j_7}$ is just an
antisymmetrization over indices $i_1,...,i_7$.

The central elements corresponding to $m=3,4$ can be expressed
through the Casimir element (since primitive central elements have
orders $2,6$) and the element corresponding to $m=7$ equals to
$detL$.

\subsection{Conclusion}

A construction of central elements in $U(\mathfrak{o}_N)$ is given.
In this construction pfaffians appear in a natural manner.

Also  a construction of central elements in $U(\mathfrak{g}_2)$ is
given. In this construction new functions of central elements $
G_{m_1,...,m_q}^{n_1,...,n_p}(L)$ given by formula \eqref{new1}
appeared.

\end{document}